\documentclass[12pt]{article}
\begin{document}

\newtheorem{prop}{Proposition}
\newtheorem{defi}[prop]{Definition} \newtheorem{thm}[prop]{Theorem} 
\newtheorem{lemma}[prop]{Lemma} \title{First 
neighbourhood of the diagonal, and geometric  distributions} 
\author{Anders Kock} \maketitle

\noindent The consideration of the $k$'th neighbourhood of the diagonal of 
a manifold $M$, 
$M_{(k)} \subseteq M\times M$,
was initiated by Grothendieck to import notions from 
differential geometry into the realm of algebraic geometry. 
These notions were re-imported into differential geometry by 
Malgrange \cite{malgrange}, Kumpera and Spencer, \cite{KS}, \ldots.  They 
utilized the notion of {\em ringed} space (a space equipped with a {\em 
structure sheaf} of functions).  The only points of (the underlying 
space of) $M_{(k)}$ are the diagonal points $(x,x)$ with $x\in M$.  
But it is worthwhile to describe mappings to and from $M_{(k)}$ {\em 
as if} $M_{(k)}$ consisted of ``pairs of $k$-neighbour points $(x,y)$" (write 
$x\sim _k y$ for such a pair; such $x$ and $y$ are ``point proches" in 
the terminology of A.  Weil).  The introduction of topos theoretic 
methods has put this ``synthetic" way of speaking onto a rigourous 
basis, and we shall freely use it.

We shall be only interested in the case $k=1$, so we are considering the 
{\em first} neighbourhood of the diagonal, $M_{(1)} \subseteq M\times M$, 
and we shall write $x\sim y$ instead of $x\sim _1 y$ whenever 
$(x,y)\in M\times M$ belongs to the subspace $M_{(1)}$. The relation $\sim$ is 
reflexive and symmetric, but unlike the notion of ``neighbour'' relation in 
Non Standard Analysis, it is {\em not} transitive.  In some previous 
writings, \cite{SDG}, \cite{DFVG}, \cite{CCBI}, \cite{DFIC}, \cite{APFBAC}, 
we have discussed the paraphrasing of several differential-geometric 
notions in terms of the combinatorics of the neighbour relation 
$\sim$.  We shall remind the reader about some of them concerning 
differential forms and connections in Section 5 below.

The new content in the present note is a description of the 
differential-geometric notion of {\em distribution} and {\em 
involutive} distribution in combinatorial terms (Section 3, notably 
Theorem \ref{involutive}), and an application of these notions to a 
proof of a version of the Ambrose-Singer holonomy theorem (Section 5).

 The proofs will depend on a correspondence between the classical Grassmann 
 algebra, and the algebra of ``combinatorial differential forms'', 
 which has been considered in the series of texts by the author, 
 mentioned above, and by Breen and Messing, \cite{BM}.  We summarize, 
 without proofs, the relevant part of this theory in the first two 
 sections.  (We hope that this summary may have also an independent 
 interest.)

\section{The algebra of combinatorial \\ differential forms}

In the smooth category, a manifold $M$ can be recovered from the ring 
$C^{\infty}(M)$ of smooth functions on it, (-- unlike in the analytic category, 
or the category of schemes 
in algebraic geometry, where there are not enough {\em global} 
functions on $M$, so that one here needs the {\em sheaf} of (algebraic) 
functions to recover the space).  In any case, ``once the seed of 
algebra is sown it grows fast'' (Wall, \cite{W} p.\ 185): if we let a 
geometric object like a manifold be represented by a ring (its ``ring 
of functions''), we may conversely consider (suitable) rings to 
represent ``virtual'' geometric objects.  The ring in question 
(better: {\bf R}-algebra, or even ``$C^{\infty}$-algebra'', cf.\ e.g.\ 
\cite{MR}) is then considered as the ring of smooth real valued 
functions on the virtual object it represents.  In particular, from 
Grothendieck, Malgrange and others, as alluded to in the introduction, 
we get the following geometric object $M_{(1)}$, where $M$ is a 
manifold:

In the ring $C^{\infty} (M\times M)$, we have the ideal $I$ of functions 
$f(x,y)$ vanishing on the diagonal, $f(x,x)=0$, in other words $I$ is the 
kernel of the restriction map $C^{\infty} (M\times M) \to C^{\infty}(M)$ 
(restrict along the diagonal $M\to M\times M$). Since this restriction map is 
a surjective ring homomorphism, we have $C^{\infty}(M) \cong  C^{\infty} (M\times 
M)/I$.  Now $I$ contains $I^2$, the ideal of functions vanishing to 
the {\em second} order on the diagonal.  The ring $C^{\infty}(M\times 
M)/I^2$ represents a new, virtual, geometric object $M_{(1)}$, and the 
surjection $C^{\infty}(M\times M)/I^2 \to C^{\infty}(M\times M)/I 
\cong C^{\infty}(M)$, we may see as restriction map for a ``diagonal'' 
inclusion $M\to M_{(1)}$.

The points of $M$ may be recovered from the algebra $C^{\infty}(M)$ as the 
algebra maps $C^{\infty}(M)\to{\bf R}$. The ``virtual'' character of $M_{(1)}$ 
is the fact that the {\em real} points of it are the same as the real points 
of $M$, or put differently, a real point $(x,y)$ in $M\times M$ belongs to 
$M_{(1)}$ precisely when $x=y$. So we cannot get a full picture of $M_{(1)}$ 
by looking at its real points. The synthetic way of speaking, and reasoning, 
talks about objects like $M_{(1)}$ in terms of the virtual points $(x,y)$,
and this mode of speaking can be fully justified by interpretation via 
categorical logic in suitable ``well adapted toposes'', see e.g.\ 
\cite{D},  \cite{SDG}. 
To illustrate the synthetic language concretely, let us consider the 
following definition.

\begin{defi}A {\em combinatorial differential 1-form} $\omega$ on a 
manifold $M$ is a function $M_{(1)} \to {\bf R}$ which vanishes on the 
diagonal $M\subseteq M_{(1)}$.
  \label{one} \end{defi}
   
 Let us analyze this notion in non-virtual terms. A {\em function}  $\omega : M_{(1)}\to {\bf 
 R}$ is  an element of the ring $C^{\infty }(M\times M)/I^2 $, 
 by the very definition. When does such $\omega$ satisfy the clause 
 that ``it vanishes on the diagonal''?  precisely when it belongs to 
 $I$ (-- the functions vanishing on the diagonal), or more precisely, 
 when it belongs to the image of $I$ under the homomorphism 
 $C^{\infty}(M\times M) \to C^{\infty}(M\times M)/I^2 $.  So a 
 combinatorial differential 1-form on $M$ is an element of $I/I^2$; 
 but this is the module of {\em Kaehler differentials} on $M$.  Kaehler 
 observed that this module was isomorphic to the classical module of 
 differential 1-forms on $M$, i.e.\ the module of fibrewise linear maps 
 $\overline{\omega}: T(M)\to {\bf R}$.  This results thus can be 
 formulated:
 
 \begin{thm}(Kaehler) There is a bijective correspondence between the mo\-dule 
 of classical 
 differential 1-forms on $M$, and the module of functions $M_{(1)} \to {\bf 
 R}$ which vanish on the diagonal. 
 Or, in the terminology of Definition \ref{one}: a 1-1 correspondence between 
 classical differential 1-forms $\overline{\omega }$ and combinatorial 
 differential 1-forms $\omega$.  \end{thm}
 
Note that for combinatorial 1-forms, there is {\em no} linearity requirement. 
Furthermore, combinatorial 1-forms are automatically {\em alternating} in the 
sense that $\omega (x,y)= -\omega (y,x)$ for all $x\sim y$. This is 
essentially because for any function $f(x,y)$ vanishing on the diagonal, the 
function $f(x,y)+ f(y,x)$ also vanishes on the diagonal and is 
furthermore symmetric in $x$, $y$, hence vanishes to the second order 
on the diagonal.

\medskip

There are also combinatorial versions of the notion of differential 
$k$-form, $k\geq 2$, \cite{SDG}, \cite{CCBI}, \cite{DFIC}, \cite{BM}, 
and a comparison result with 
classical (multilinear alternating) differential forms. We summarize the 
definition and comparison here, using synthetic language. A $k+1$-tuple of 
elements $(x_0 , \ldots ,x_k )$ of (virtual !) points of $M$ is called an {\em 
infinitesimal $k$-simplex} if $x_i \sim x_j$ for all $i$, $j$, and is called 
{\em degenerate} if $x_i = x_j$ for some $i\neq j$. The ``set'' of 
infinitesimal $k$-simplices in $M$ form an object $$M_{[k]}\subseteq M^{k+1}$$ 
(which may be described in ring theoretic terms, in the spirit of $I/I^2$, but 
more complicated; cf.\  \cite{BM}).
The following definition from \cite{SDG} is then an extension of the previous 
definition:

\begin{defi} A combinatorial differential $k$-form  $\omega$ on a manifold $M$ is 
a function $M_{[k]} \to {\bf R}$ which vanishes on all degenerate $k$-simplices.
   \end{defi}
   
   In the rest of this text, we shall often just say (combinatorial) $k$-form instead of 
   combinatorial differential $k$-form.
   
Note again that in the Definition, there is no (multi-)linearity requirement; and $\omega$ can 
be proved to be alternating in the sense that the interchange of $x_i$ and 
$x_j$ results in a change of sign in the value of $\omega$. 

We proceed to describe exterior derivative of combinatorial forms, and wedge 
product; both these structures are analogous to structures (coboundary and cup 
product) on the singular cochain complex of a 
topological space (see \cite{DFIC}). 

For $\omega$ a $k$-form, we let $d\omega$ be the $k+1$-form given by
$$d\omega \; (x_0 , \ldots ,x_{k+1}):= \sum _{i=0} ^{k+1} (-1)^i \omega 
(x_0 ,\ldots ,\hat{x_i}, \ldots ,x_{k+1});$$
for $\omega$  a $k$-form and for $\theta$ an $l$-form, we let $\omega 
\wedge \theta$ be the $k+l$-form given by
$$(\omega 
\wedge \theta )(x_0 , \ldots , x_{k+l}):= \omega (x_0 , \ldots ,x_k )\cdot 
\theta (x_k ,x_{k+1}, \ldots ,x_{k+l}).$$

(It is not trivial, but true, that one gets the value 0 when $\omega \wedge 
\theta$ is applied to a $k+l$-simplex which is degenerate by virtue of $x_i = 
x_j$ with $i< k$ and $j>k$.)

Equipped with these structures $d$ and $\wedge$, together with the evident 
``pointwise'' vector space structure, the combinatorial forms together make up 
a differential graded algebra $\Omega ^{\bullet} (M)$. Let 
$\overline{\Omega}^{\bullet}(M)$ be the classical differential graded algebra 
 of classical differential forms. Part of the content of \cite{SDG}, 
\cite{CCBI}, \cite{DFIC} may be summarized in the following ``Comparison'' 
Theorem (see also \cite{BM} for the generalization to  schemes):

\begin{thm}The differential graded algebras $\Omega ^{\bullet} (M)$ and 
$\overline{\Omega}^{\bullet}(M)$ are isomorphic.   
  \label{main} \end{thm}
   
   The correspondence is given in \cite{SDG}, Corollary 18.2; the 
   compatibility with the differential and product structure is proved in 
   \cite{DFIC} p.\ 259 and 263, resp.
 
It should be said, though, that the strict validity of the theorem depends on 
the conventions for defining the classical exterior derivative and wedge; the 
conventions often differ by a factor $k!$  or $k! l! /(k+l)!$. So with 
different conventions, the equalities claimed by the theorem hold only modulo 
such rational factors. The uses we shall make of the theorem, however, have 
the character ``if one thing is zero, then so is the other'', 
so are independent on these rational factors, and thus indpendent of 
the conventions.

\medskip
\section{The $D$-construction; $\log$ and $\exp$}

For any finite dimensional vector space $E$, we let $D(E)$ denote the 
subset of $x\in E$ with $x\sim 0$.  If in particular $E={\bf R}$, we write 
$D$ for $D({\bf R})$.  This is the most basic object in Synthetic 
Differential Geometry, and appears as such in \cite{Lawv}, 
\cite{SAFD}.  In algebraic guise, when $D$ is presented in terms of 
the ring of functions on it, it appears much earlier, namely as the 
``ring ${\bf R}[\epsilon ] = {\bf R}[X]/(X^2  )$ of dual numbers''.  The object 
$D\subseteq {\bf R}$ is the set of (virtual) elements $d\in {\bf R}$ 
which satisfy $d^2 =0$.  If $M$ is a manifold, a map $t: D\to M$ with 
$t(0)=x \in M$ is a {\em tangent vector} at $x\in M$.  (This assertion 
is just the synthetic reformulation of: an algebra map 
$C^{\infty}(M)\to {\bf R}[\epsilon]$ corresponds to a derivation 
$C^{\infty}(M)\to {\bf R}$.)

\medskip

A fundamental fact about ``ringed spaces'' gets its synthetic 
formulation in the following  ``cancellation principle'':

\begin{equation}\mbox{ If $a\in {\bf R}$ satisfies $d\cdot a =0$ for {\em 
all} $d\in D$, then $a=0$.  }\label{cancellation}\end{equation}

(Verbally: {\em universally quantified $d$'s may be cancelled''}. It is 
part of what sometimes goes under the name ``Kock-Lawvere axiom'', see 
e.g.\  \cite{MR}.)

From this cancellation principle, it is immediate to deduce the 
following

\begin{prop} Let $E$ be a vector space, and let $U$ and $V$ be linear 
subspaces.  If $V$ is the kernel of a linear map $\phi : E \to {\bf 
R}^k$, and if $d\cdot U \subseteq V$ for all $d\in D$, them 
$U\subseteq V$.\label{linearsufficiency}\end{prop}

Note that not every linear subspace $V$ is a kernel; for instance, the 
``set'' $D_{\infty} \subseteq {\bf R}$ of nilpotent elements is a 
linear subspace of a 1-dimensional vector space; it is not a kernel, 
and it does not qualify as finite dimensional.)

\medskip

We consider again a general finite dimensional vector space $E$. 
Let $\tilde{D}^p (E)  \subseteq 
(D (E))^p $ denote the set of $p$-tuples $( x_1 , \ldots ,x_p )$ with $x_i 
\in D(E)$ and which sa\-tisfy $x_i \sim x_j$ for all $i,j = 1, \ldots ,p$.  
So we have inclusions \begin{equation}\tilde{D}^p (E) \subseteq D(E)^p 
\subseteq E^p .\label{inclu}\end{equation} Let us for short call a 
(partial) ${\bf R}$-valued function $\theta $ on $E$ {\em normalized} 
if its value is 0 as soon as one of the input arguments is $0\in E$.  
The following auxiliary result is proved (for $E={\bf R}^n$) in 
\cite{SDG} I.16, for the bijection between 1) and 3); the bijection 
between 1) and 2) is proved in a similar way, but is easier.

 \begin{prop}Restriction along the 
inclusions in (\ref{inclu}) establish bijections between
\begin{enumerate}
\item multilinear alternating maps $E^p \to {\bf R}$ 
\item normalized 
alternating maps $D(E)^p \to {\bf R}$ 
\item normalized maps 
$\tilde{D}^p (E)  \to {\bf R}$
\end{enumerate}
\label{Extension}\end{prop}

\medskip
Now consider a manifold $M$, then its tangent bundle $TM$ is likewise 
a manifold, and we may talk about when two tangents are neighbours.  
In particular, we may ask when a tangent vector at $x\in M$ is 
neighbour to the zero tangent vector at $x$.  The set of all such 
tangents is denoted $DM \subseteq TM$.  Alternatively $DM$ is formed by 
applying the $D$ construction for vector spaces, as given above, to each 
fibre $T_x M$ individually.

\medskip

From \cite{DK}, we have that there is a canonical bijection $\exp: DM \to 
M_{(1)}$.  Its inverse, which we of course have to call $\log$, may be 
described explicitly as follows: if $(x,y)\in M_{(1)}$, then $\log 
(x,y)\in TM$ is  the tangent vector 
at $x\in M$ in $M$ given by
\begin{equation}d\mapsto (1-d)x + dy \mbox{ for } 
d\in D.\label{aff}\end{equation}
 Here we are 
forming an {\em affine combination} $(1-d)x + dy$ of two points $x$ and $y$ 
in a manifold $M$; this can be done (for {\em any} $d\in {\bf R}$, in fact)  by 
choosing a coordinate chart around $x$ and $y$, and making the affine 
combination in coordinates.  It is proved in \cite{GCLCC} Theorem 1 
that when $x\sim y$, the result does not depend on the coordinate 
chart chosen.

\medskip

For $s\in {\bf R}$ and $t\in TM$, $s\cdot t$ denotes the tangent 
vector given by $\delta \mapsto t(s\cdot \delta )$.  For $d\in D$ and 
$t\in TM$, $d\cdot t \in DM$, and also $t(d)\sim t(0)$ ($=x$, say).  We 
have the equations
\begin{equation}\exp (d\cdot t) = t(d) \mbox{ ; and } \log (x, t(d)) = d\cdot t
\label{general}\end{equation}

\medskip

In terms of $\log$, we can be more explicit about the correspondence of the 
 Theorem \ref{main}: the combinatorial $p$-form $\theta$ 
corresponding to a classical $p$-form $\overline{\theta}$ is given by 
\begin{equation}\theta (x_0 , \ldots ,x_p ) = \overline{\theta }(\log 
(x_0 , x_1 ), \ldots ,\log (x_0 , x_p )).\label{logg}\end{equation}

Note that the left hand side vanishes if $x_i = x_0$, because 
$\overline{\theta}$ is normalized (being multilinear), and vanishes if $x_i = 
x_j$ ($i, j \geq 1$) because $\overline{\theta}$ is alternating.

Note also that the right hand side is defined on more $p+1$-tuples 
than the left hand side, since on the right hand side, no assumption 
$x_i \sim x_j$ for $i ,j \geq 1$ enters.

\section{Geometric distributions}

\begin{defi} A {\em predistribution} on a manifold $M$ is a reflexive 
symmetric refinement $\approx$ of the relation $\sim$. 
\end{defi}
(That $\approx$ is a {\em refinement} of $\sim$ is taken in the sense: 
$x\approx y$ implies $x\sim y$, for all $x$ and $y$.)  For instance, if 
$f: M\to N$ is a submersion, we get a predistribution $\approx$ by 
putting $x\approx y$ when $x\sim y$ and $f(x)=f(y)$.  A 
predistribution arising in this way is clearly {\em involutive} in the 
following sense:

\begin{defi}A predistribution $\approx$ on $M$ is called (combinatorially)
   {\em involutive} if 
 $x\approx y$, $x\approx z$, and $y\sim z$ implies $y\approx z$.
\label{involutive}\end{defi}

In analogy with the object $M_{[k]}$ of infinitesimal $k$-simplices, we may in 
the presence of a distribution define $M_{[[k]]}\subseteq M_{[k]}$ to consist 
of $k+1$-tuples $(x_0 , \ldots ,x_k )$ satisfying $x_i \approx x_j$ for all 
$i$ and $j$. If we call such infinitesimal simplices {\em flat}, we can 
reformulate the notion of an involutive distribution as follows: {\em if two 
of the faces of an infinitesimal 2-simplex are flat, then so is the third}; 
or equivalently: {\em if two 
of the faces of an infinitesimal 2-simplex are flat, then so is the 2-simplex 
itself}.

\medskip

We need to explain in which sense a (geometric) distribution on $M$ in the 
classical sense gives rise to a predistribution.  Recall that a 
$k$-dimensional distribution on an $n$-dimensional manifold $M$ is a 
$k$-dimensional sub-bundle $E$ of the tangent bundle $TM\to M$.  For 
$x\in M$, $E_x \subseteq T_x M$ is thus a $k$-dimensional linear 
subspace of the tangent vector space to $M$ at $x$.

 \medskip

If $E \subseteq TM$ is a distribution in the classical sense, we can 
for $x\sim y$ define a predistribution $\approx$ (or  $ \approx _E$) by 
\begin{equation}x\approx y \mbox{ iff }\log(x,y) \in E_x 
,\label{approx}\end{equation} 
or equivalently, the set of $y$ ($\sim 
x$) which satisfy $x\approx y$ is the image of $E_x \cap DM$ under 
$\exp$.  It is clearly a reflexive relation, since $\log (x,x)$ is the 
zero tangent vector at $x$.  The symmetry of $\approx$ is less 
evident, but is proved in Proposition \ref{symm} below.  A 
predistribution which comes about in this way from a classical 
distribution $E$, we call simply a {\em distribution}.  (This is 
justified, since one can prove (using Proposition \ref{linearsufficiency}) 
that two classical distributions, giving rise to the same 
predistribution, must themselves be equal. 

\begin{prop} Let $\approx$ be derived from a classical distribution, as in 
(\ref{approx}). Then $x\approx y$ implies $y\approx x$.
\label{symm} \end{prop}

{\bf Proof.} Recall that a classical distribution $E$, of dimension $k$, say, may 
be presented locally by $n-k$ non-singular differential 1-forms on $M$,  
$\overline{\omega} _1, \ldots , \overline{\omega} _{n-k}$, with $t\in 
E_x$ iff $t$ is annihilated by all the $\overline{\omega} _i$'s. So $x\approx 
y$, iff $\overline{\omega }_i (\log (x,y)) =0$ for $i= 1,\ldots ,(n-k)$, 
iff $\omega _i (x,y) =0$ for $i= 1,\ldots ,(n-k)$. But since each $\omega _i$ 
is alternating, the relation $\omega _i (x,y) =0$ is symmetric.

\medskip

   The justification for our use of the phrase ``involutive'' is contained 
in

\begin{thm} Let the classical distribution $E\subseteq TM$ be given.  
Then it is involutive in the classical sense if and only if  $\approx _E$ is 
involutive in the combinatorial sense of Definition \ref{involutive}.
\end{thm}

{\bf Proof.} Present, as above, the classical distribution $E$ by $n-k$ 
differential 1-forms on $M$ 
$\overline{\omega} _1, \ldots , \overline{\omega} _{n-k}$.  If 
$\overline{I}$ denotes the ideal in the Grassman (exterior) algebra 
$\overline{\Omega}^{\bullet}(M)$ generated by 
the $\overline{\omega} _i$'s, then by the classical definition, the 
distribution is involutive precisely when $\overline{I}$ is closed 
under exterior differentiation.  (It suffices that each 
$d\overline{\omega} _i$ is in $I$.) (The alternative, equivalent, 
definition of involutiveness in terms of vector fields along $E$ is 
less convenient for the comparison we are about to make.) 
 
 Now assume that $E$ is classically involutive.  Let $x\approx y$, 
 $x\approx z$ and $y\sim z$.  We need to prove $y\approx z$.  It 
 suffices to prove for each $i$ that $\overline{\omega}_i (\log (y,z)) 
 =0$, or equivalently that $\omega _i (y,z)=0$, where the $\omega _i$ 
 corresponds to $\overline{\omega}_i$ under the correspondence of 
 Theorem \ref{main}.  But by assumption $d \overline{\omega}_i \in 
 \overline{I}$, hence $d\omega _i \in I$, and so
 $$0=d\omega _i (x,y,z) = \omega _i (x,y) + \omega _i (x,z ) - \omega _i (y,z).$$
 Since the two first terms here are $0$ by $x\approx y$ and $x\approx 
 z$, we conclude $\omega _i (y,z) =0$.  So $\approx$ is 
 combinatorially involutive.
 
Conversely, assume that $\approx _E$ is combinatorially involutive.  
Let $\overline{\omega}$ be one of the 1-forms generating 
$\overline{I}$. Being 
generated by 1-forms, it follows from classical multilinear algebra 
that the ideal $\overline{I}$ has the property: for $\overline{\theta}$  
a $p$-form, if
$$\overline{\theta} : TM \times _M \ldots \times _M TM \to {\bf R}$$
 annihilates $E^p$, then $\overline{\theta} \in \overline{I}$.  So to 
 prove $d\overline{\omega} \in \overline{I}$, it suffices to prove that 
 if $t_1 , t_2 \in E_x \subseteq T_x M$, then
$$d\overline{\omega} (t_1 , t_2 ).$$
By bilinearity, and by the cancellation principle (\ref{cancellation}) 
applied twice, it suffices to prove, for all $d_1$ and $d_2$ in $D$, that 
$0=d\overline{\omega} (d_1\cdot t_1 , d_2 \cdot t_2 ).$ We calculate 
this expression, using (\ref{general}): $$d\overline{\omega} (d_1\cdot 
t_1 , d_2 \cdot t_2 ) = d\overline{\omega} (\log (x, t_1 (d_1 )), \log 
(x, t_2 (d_2 ))$$ $$ = \overline{d\omega} (\log (x, t_1 (d_1 )), \log (x, 
t_2
(d_2 ))$$
where $d\omega$ is the combinatorial exterior derivative of $\omega$, and 
$\overline{d\omega}$ the corresponding classical form (using that the 
bijection $\theta \leftrightarrow  \overline{\theta}$ of Theorem \ref{main} 
commutes with exterior derivative).  Unfortunately, we do not 
necessarily have $t_1 (d_1 ) \sim t_2 (d_2 )$, so that we cannot 
rewrite this as $d\omega (x, t_1 (d_1 ), t_2 (d_2 ))$; we need first 
an auxiliary consideration for the 2-form $\theta = d\omega $.
 
 \medskip

For $\theta$ a combinatorial 2-form, we consider a certain  function 
$\tilde{\theta}$ defined on all ``semi-infinitesimal 2-simplices'', 
meaning triples $x,y,z$ with $x\sim y$ and $x\sim z$ (but not 
necessarily $y\sim z$). The function $\tilde{\theta}$ is defined by
$$\tilde{\theta} (x,y,z) := \overline{\theta}(\log (x,y),\log (x,z)).$$
If $y\sim z$, this value is just $\theta (x,y,z)$, so 
$\tilde{\theta}$ is an extension of $\theta$.

\begin{lemma}If $\theta$ annihilates all infinitesimal 2-simplices 
$(x,y,z)$ with $x\approx y$ and $x\approx z$, then $\tilde{\theta}$
annihilates all semi-infinitesimal 2-simplices $(x,y,z)$ with $x\approx 
y$ and $x\approx z$.
\end{lemma}

{\bf Proof.} Let the semi-infinitesimal 2-simplex $x,y,z$ be given.  
Let $V$ denote $T_x M$ and let $E$ denote $E_x \subseteq T_x M$. Under the 
$\log /\exp $-identification, we identify the set of $\{y \mid  y\sim x\}$ with 
$D(V)$, and the set of $\{y \mid y\approx x\}$ with $D(E)$.  Then $\theta (x, 
-, -)$ and $\tilde{\theta }(x,-,-)$ are functions $\theta : 
\tilde{D}^2 (V) \to R$ and $\tilde{\theta}: D(V)^2 \to V$.  The 
function $\theta$ is normalized, $\tilde{\theta}$ is normalized and 
alternating.  We have similarly the restrictions of $\theta$ and 
$\tilde{\theta}$ to $\tilde{D}^2 (E)$ and $D(E)^2$.  We may see 
$\tilde{\theta}: D(E)^2 \to {\bf R}$ as arising {\em either} by first 
extending $\theta : \tilde{D}^2 (V) \to {\bf R}$ to $D(V)^2 \to {\bf 
R}$, and then restricting it to $D(E)^2$, {\em or} as arising by first 
restricting to $\tilde{D}^2 (E)$ and then extending.  Using the 
correspondence between items 2) and 3) in Proposition 
\ref{Extension}, we conclude that this must give the same result.  
The assumptions made on $\theta$ guarantees that the ``first 
restricting, then extending'' process yields 0, hence so does the 
other process, so $\tilde{\theta}$ restricts to 0 on $D(E)^2$; but 
this is, under the $\log /\exp$ identification, precisely the 
assertion that $\tilde{\theta}(x,y,z) =0$ for all $y\approx x$ and 
$z\approx x$.

\medskip

We can now finish the proof.  We have that the 2-form $d\omega$ 
annihilates infinitesimal 2-simplices $(x,y,z)$ with $x\approx y$ and 
$x\approx z$, because then, by assumption, $y\approx z$. From the Lemma, it then 
follows that $\widetilde{d\omega}$ annihilates semi-infinitesimal 
2-simplices $x,y,z$ with with $x\approx y$ and $x\approx z$.  Since 
$t_1 (d_1 )\approx x$ and $t_2 (d_2 ) \approx x$, the simplex $(x, t_1 
(d_1 ), t_2 (d_2 ) )$ is of this kind, so $\widetilde{ d\omega}$ takes 
value 0 on this simplex.  But $$\widetilde{d\omega}(x, t_1 (d_1 ), t_2 
(d_2 )) = \overline{d\omega} (\log (x, t_1 (d_1 ), \log (x, t_2
(d_2 )),$$
which thus is 0.  This proves that classical involution condition.

\section{Frobenius Theorem}
Let $E\subseteq TM$ be a distribution on a 
manifold $M$.

\begin{prop}Let $E\subseteq TM$ be a classical distribution, and $\approx$ 
the corresponding combinatorial one. Let $Q\subseteq M$ be a 
submanifold, and $x\in Q$
Then the following are equivalent: 

1) for any $t\in T_x M$: $ t\in T_x Q$ implies $t\in E_x$

2) for any $y$ with $x\sim y$, $y\in Q$ implies $x\approx y$.

\noindent Also the following are equivalent:

3) for any $t\in T_x M$, $t\in E_x$ implies $t\in T_x Q$

4) for any $y$ with $x\sim y$,  $x\approx y$ implies $y\in Q$.

\label{two}\end{prop}

{\bf Proof.} Assume 1). Let $x\sim y$, $y\in Q$. Then $\log 
(x,y) \in T_x Q$ (as a submanifold, $Q$ is stable under the 
formation of the affine combinations that make up the values 
of $\log$), so $\log (x,y) \in E_x$, meaning $x\approx y$.

Conversely, assume 2). To prove $T_x Q \subseteq E_x$, it 
suffices by Proposition \ref{linearsufficiency}  to prove for each $d\in D$ that 
$d\cdot T_x Q \subseteq E_x$.  So consider a tangent vector $d\cdot t$ 
where $t\in T_x Q$.  Then $t(d)\in Q$ and $x\sim t(d)$, whence by 
assumption $x\approx t(d)$, i.e.\ $\log (x,t(d)) \in E_x$.  But $\log 
(x, t(d))=d\cdot t$ by (\ref{general}), so $d\cdot t \in E_x$.

Assume 3), and assume $x\approx y$, i.e.\ $\log (x,y)\in E_x$. 
By assumption, therefore,  $\log (x,y)\in T_x Q$. Applying 
$\exp$ yields $y\in Q$.

Conversely, assume 4). To prove $E_x \subseteq T_x Q$,
 let $t\in E_x$.  Then by (\ref{general}), for each $d\in D$, $x\approx \exp 
 ( d\cdot t)$, whence by assumption $\exp (d\cdot t)\in Q$, so $t(d)\in 
 Q$, by (\ref{general}).  Since this holds for all $d\in D$, $t$ is a 
 tangent vector of $Q$.
 
\medskip 

 Consider a classical distribution $E\subseteq TM$, and let 
$\approx$ be the corresponding combinatorial one.  A submanifold 
$Q\subseteq M$ is called {\em weakly integral} if for each $x\in Q$, 
the equivalent conditions 1) and 2) of Proposition \ref{two} holds; 
and is called (strongly) {\em integral} if furthermore, the equivalent 
conditions 3) and 4) hold.

\medskip

The Frobenius Theorem says that if $E\subseteq TM$ is an 
involutive distribution, then there exists for every $x\in M$ 
a (stronly) {\em 
integral} submanifold $Q\subseteq M$ containing $x$; there even exists 
a unique maximal connected such $Q$.  Maximality here means that if 
$x\in K\subseteq M$ is any connected submanifold which is weakly integral, then 
$K\subseteq Q$. 
 \medskip
 
 We shall later use the following ``sufficiency'' principle \begin{prop}Let 
 $H\subseteq G$ be a Lie subgroup of a connected Lie group.  Let ${\cal 
 M}(e)$ be the set of $g\in G$ with $g\sim e$ ($e$ the neutral 
 element).  If ${\cal M}(e)\subseteq H$, we have $G=H$.
 \label{sufficiency}\end{prop}

{\bf Proof.} Finite dimensional linear subspaces of a finite 
dimensional vector space $E$ may be reconstructed from their meet with 
$D(E)$, by Proposition \ref{linearsufficiency}.  Using the $\log$ - $\exp $ 
bijection, finite dimensional linear subspace of $T_e (G)$ may be 
reconstructed from subsets of ${\cal M}(e)$.  The assumption of the 
Proposition therefore gives that $T_e H = T_e G$, and this implies by 
classical Lie theory that $H=G$.

\section {Ambrose-Singer Theorem} This Theorem 
(\cite{AS}, or see \cite{Nomizu}, II.7), deals with connections in principal 
bundles.  We briefly recall how this gets formulated in synthetic 
terms; for a more elaborate account, see \cite{CCBI} and \cite{APFBAC}.  
Let $\pi : P\to M$ be a principal $G$-bundle, in the smooth category 
of course; $G$ is a Lie group.  The action of $G$ on $P$ is from the 
right.  A {\em principal connection} $\nabla$ in $P\to M$ is a law 
which to any $a\sim b$ in $M$ assigns a $G$-equivariant map $\nabla 
(a,b): P_b \to P_a$, with $\nabla (a,a)$ the identity map (this 
implies that $\nabla (a,b)$ has $\nabla (b,a)$ as an inverse).  The 
{\em connection form} $\omega$ of $\nabla$ is a ``$G$-valued 1-form on 
$P$'': it is the law which to any $x\sim y \in P$ assigns an element 
of $G$, namely the element $\omega (x,y)\in G$ such that
$$x\cdot \omega (x,y) = \nabla (a,b)(y),$$
where $a =\pi (x)$ and $b=\pi (y)$.  If one thinks of an 
infinitesimal 1-simplex of the shape $( \nabla (a,b)(y), y)$ as a {\em 
horizontal} 1-simplex, $\omega (x,y)$ thus measures the lack of 
horizontality of $(x,y)$; in particular, $(x,y)$ is horizontal 
precisely when $\omega (x,y)=e$, the neutral elemnt of $G$.  A curve 
in $P$ is called horizontal if any two neighbour points on it form a 
horizontal 1-simplex.

 The {\em 
curvature form} $d\omega$ is the $G$-valued exterior derivative of 
$\omega$, namely the law which to an infinitesimal 2-simplex $(x,y,z)$ 
in $P$ assigns the element $\in G$ 
\begin{equation}d\omega (x,y,z):= 
\omega (x,y)\cdot \omega (y,z) \cdot \omega 
(z,x).\label{coboundary}\end{equation}

\medskip

The following is now a (restricted) version of the Ambrose Singer ``holo\-nomy'' 
theorem.  A principal $G$-bundle with connection is given, as above; 
$G$ is assumed connected.
 
\begin{thm} Assume that any two points of $P$ can be connected by a 
horizontal curve, and that $H\subseteq G$ is a Lie subgroup such that 
for all infinitesimal 2-simplices $(x,y,z)\in P$, we have $d\omega 
(x,y,z)\in H$.  Then $H=G$.
\end{thm}

{\bf Proof.} We construct a distribution $\approx$ on $P$, by putting 
$$x\approx y \mbox{ iff } \omega (x,y )\in H.$$ 
In particular, if $x$ and $y$ are in the same fibre $P_a$,
\begin{equation}x\approx y \mbox{ iff } y=x\cdot 
h\label{pentagram}\end{equation}
 for some $h\in H$ with  $h\sim e$, ($e$ the neutral element of $G$).

 This distribution 
is clearly involutive; for if $(x,y,z)$ is an infinitesimal 2-simplex 
with $x\approx y$ and $x\approx z$, three of the 
four factors that occur in the equation (\ref{coboundary}) are in $H$, 
hence so is the fourth factor $\omega (y,z)$, proving $y\approx z$.

Pick a point $x\in P$, with $\pi (x)=a$, say.  By Frobenius, there is 
a maximal integral submanifold $Q\subseteq P$ for the distribution 
$\approx$.  Now horizontal infinitesimal 1-simplices $(x,y)$ have 
$x\approx y$, since $\omega (x,y)=e$ (the neutral element of $G$).  
Therefore by the Proposition \ref{two} and maximality of $Q$, any 
horizontal curve through $x$ must lie entirely in $Q$. Since any 
point of $P$ can be connected to $x$ by a horizontal curve, $Q=P$.

Let ${\cal M}(x)$ denote the set of neighbours of $x$.  Since $Q=P$, 
we then have 
$$P_a \cap {\cal M}(x) = Q_a \cap {\cal M}(x) =\{y\in P_a \mid 
y\approx x\}$$ since $Q$ is a integral manifold for $\approx$.  But this 
 is $\{ y\sim x \mid y=x\cdot h \mbox{ for some } h\in H$, by 
({\ref{pentagram}).  The left hand side here consists of elements 
$x\cdot g$ with $g\sim e$.  So $g\sim e$ in $G$ implies $g\in H$.  
From Proposition \ref{sufficiency}, we deduce that $H=G$.


\begin{thebibliography}{99}

\bibitem{AS} W.\ Ambrose and I.\ Singer, A theorem on 
holonomy, Trans.\ Amer.\ Math.\ Soc.\ 75 (1953), 428-443.
\bibitem{BM}L.\ Breen and W.\  Messing, Combinatorial differential forms, 
Adv.\ Math.\ 164  (2001), 203-282.
 \bibitem{D}E.\ Dubuc, Sur les modeles de la geometrie 
 differentielle synthetique, Cahiers de Top.\ et Geom.\ diff.\ 
 20 (1979), 231-279.
\bibitem{DK}E.\ Dubuc and A.\ Kock, On 1-form classifiers, 
Comm.\ in Algebra 12 (1984), 1471-1531.
\bibitem{SAFD} A.\ Kock, A simple axiomatics for differentiation, 
Math.\ Scand.\  40 (1977), 183-193. \bibitem{SDG}A.\ Kock, Synthetic Differential Geometry, 
London Math.\ Soc.\ Lecture Notes Series 51, Cambridge Univ.\ Press 
1981.  \bibitem{DFVG}A.\ Kock, Differential forms with values in 
groups, Bull.\ Austral.\ Math.\ Soc.\ 25 (1982), 357-386.
\bibitem{CCBI} 
A.\ Kock, Combinatorics of curvature and the Bianchi identity, 
Theory and Applications of Categories 2 No.\ 7 (1996).  
  \bibitem{GCLCC}A.\ Kock, Geometric constructon of the 
Levi-Civita parallelism,    Theory and Applications of 
Categories 4 No.\ 9 (1998).
\bibitem{DFIC} A.\ Kock, Differential forms as infinitesimal cochains, J.\ Pure 
Appl.\ Alg.\ 154 (2000), 257-264.  \bibitem{APFBAC}A.\ Kock, Algebra 
of principal fibre bundles, and connections, arXiv:math.CT/0005125, 
May 2000.  \bibitem{KS} A.\ Kumpera and D.\ Spencer, Lie Equations 
Volume I: General Theory, Annals of Math.\ Studies no.\ 73, 1972.  
\bibitem{Lawv} F.W.\ Lawvere, Categorical dynamics, in Topos Theoretic 
Methods in Geometry, ed.\ A.\ Kock, Aarhus Various Publications Series 
No.\ 30 (1979).  \bibitem{malgrange}B.\ Malgrange, Equations de Lie. I, 
J.\ Diff.\ Geom.\ 6 (1972), 503-522. \bibitem {MR} I.\ Moerdijk and G.E.\ Reyes, Models 
for Smooth Infinitesimal Analysis, Springer 1991.  \bibitem{Nomizu}K.\ 
Nomizu, Lie groups and differential geometry, The Mathematical Society 
of Japan 1956. \bibitem{W}C.T.C.\ Wall, Lectures on 
$C^{\infty}$-Stability and Classification, Proceedings of Liverpool 
Singularities Symposium I, 178-206, Springer Lecture Notes in Math.\ 
192 (1971).
\end{thebibliography}
\end{document}